\documentclass[12pt]{article}
\usepackage{amsthm,amsfonts,amsmath,amssymb}
\usepackage{times}

\baselinestretch

\def\EE{{\mathbb E}}

\def\RRR{{\mathbb R}}
\def\SSS{{\mathbb S}}

\title{~~~\\[-3cm]
Rate of convergence of random polarizations}
\author{Almut Burchard\footnote{
{\tt almut@math.toronto.edu}}\\
University of Toronto
}
\date{August 23, 2011}

\begin{document}

\maketitle
\thispagestyle{empty}

\begin{abstract} After $n$ random polarizations
on $\SSS^d$, the expected symmetric difference of a Borel set from a
polar cap is bounded by $C_dn^{-1}$, where
the constant $C_d$ depends on the dimension~\cite{BF}. 
We show here that this power law is best possible 
and that necessarily $C_d\ge d$.
\end{abstract}

\section{Introduction}

Let $A$ be a subset of the $d$-dimensional sphere $\SSS^d$
(viewed as the unit sphere in $\RRR^{d+1}$),
and 
let $\sigma:x\mapsto \bar x$
be a reflection at a great circle that does
not pass through the north pole.
The {\bf polarization} of $A$ 
with respect to $\sigma$ is defined by
$$
x\in SA \ \Leftrightarrow\ 
\left\{\begin{array}{ll}
   x\in A\ \ \mbox{or}\ \ \bar x\in A\,,  & \text{if } 
\delta(x,O) \le \delta(\bar x,O)\,, \\
   x\in A \ \mbox{and}\ \bar x \in A\,, & \text{if } 
\delta(x,O) \ge \delta(\bar x,O) \,.
\end{array}\right.
$$
Here, $ \delta(x,y)$ denotes the geodesic distance on $\SSS^d$ 
given by the angle enclosed between $x$ and $y$,
and $O$ denotes the north pole. Since reflections preserve
the uniform probability measure $m(\cdot)$ on
the sphere, so do polarizations, and 
\begin{equation}
\label{eq:pol-identity}
m(SA\cap SB)-m(A\cap B)  = \int_{\SSS^d}
I_{A\setminus B}(x) I_{B\setminus A}(\bar x)\,dx \ge 0\,.
\end{equation}
We parametrize the reflections on $\SSS^d$ by 
$u\in\Omega=\SSS^d/\pm$, setting 
$$
\sigma_u(x)=x-2(u\cdot x)u\,,
$$
and we denote the corresponding polarization 
by $S_u$. A {\bf random polarization} $S_U$
is polarization in the direction of a uniformly distributed 
random variable $U$ on~$\Omega$. 
We consider sequences of random polarizations 
$S_{U_1\dots U_n}=S_{U_n}\circ \cdots \circ S_{U_1}$,
where the $\{U_i\}_{i\ge 1}$ are independent.  
Van Schaftingen has shown that almost surely, 
for every Borel set $A$ the sequence $S_{U_1\dots U_n}A$ 
converges to the polar cap $A^*$ of the same 
volume~\cite[Theorem 3.13]{VS}. 
The convergence occurs in symmetric difference if $A$ is measurable, and
in Hausdorff metric if $A$ is compact~\cite[Corollary 4]{BF}.  

Subject of this note is the rate of convergence.
In prior work, we have shown that under a 
similar sequence of random polarizations on $\RRR^d$,
\begin{equation} \label{eq:converge}
\EE\Bigl[m(S_{U_1\dots U_n} A\bigtriangleup A^*)\bigr]
\le C_d n^{-1}\,.
\end{equation}
There, $A$ is a Borel measurable subset of the unit ball, $m$ 
is Lebesgue measure (normalized so that the unit ball
has measure one),
and $C_d=d\,2^{d+1}$~\cite[Proposition 4.1]{BF}.
This rate of convergence 
is much slower than what is known for other symmetrizations.
Klartag has proved that a sequence of 3d carefully chosen 
Steiner symmetrizations in $\RRR^d$ followed by a random sequence
where each step consists of $d$ orthogonal Steiner symmetrizations
converges faster than every polynomial. The leading constant 
depends only on the dimension and
grows at most polynomially~\cite[Theorem 1.5]{K}. 
Although Klartag's result applies only to 
convex bodies, it raises the question whether the power law 
in Eq.~(\ref{eq:converge})
can be improved. For random polarizations,
the answer is negative:

\bigskip\noindent{\bf Proposition.} \ {\em For random polarizations
of a Borel set $A\subset \SSS^d$, Eq.~(\ref{eq:converge}) holds 
with $C_d=2^d$.  The $n^{-1}$ power law is optimal, and the 
sharp constant satisfies $C_d\ge d$.
}

\bigskip The proof of the proposition has two parts.
Eq.~(\ref{eq:converge}) and the upper bound on $C_d$ 
are obtained by simply adjusting
Proposition~4.1 of~\cite{BF} to the sphere.
For the lower bound on $C_d$ and to prove the
optimality of the power law, we analyze 
how spherical caps move under
polarization. If $A$ is a hemisphere,
we compare the difference of $S_{U_1\dots U_nA}$
from $A^*$ with with the order statistics of 
the uniform distribution, and relate the limiting distribution 
of $n\cdot m(S_{U_1\dots U_n}A\bigtriangleup A)$
to a Gamma distribution.  We work on the sphere rather 
than on $\RRR^d$, because the additional
symmetry simplifies the calculations.  It will be 
clear from the proofs that similar bounds hold 
for the polarization of balls in~$\RRR^d$.
Other questions remain open: {\em How quickly do the
sharp constants grow with  the dimension?  What is the 
impact of the distribution of~$U$? Can one speed up 
the convergence by adapting the sequence to~$A$?}

\bigskip\noindent{\bf Acknowledgments.} \ I wish to thank
Vitali Milman for his comments on~\cite{BF}
that prompted this work, and Marc Fortier for patience and
helpful discussions.  The research was supported in part by 
NSERC through Discovery Grant No. 311685-10.

\section{The upper bound on the sharp constant}

In this section, we show that Eq.~(\ref{eq:converge})
holds on the sphere  with $C_d=2^d$.
For a single random polarization we have
by the identity (\ref{eq:pol-identity})
and Fubini's theorem,
$$
m\bigl( A\bigtriangleup A^*\bigr) - 
\EE\bigl[ m\bigl(S_U A\bigtriangleup A^*\bigr)\bigr]
= 2 \int_{A^*\setminus A} 
P(\sigma_U(x)\in A\setminus A^*)\, dm(x)\,.
$$
We compute the probability under the integral
as an average over the hemisphere where $u\cdot x >0$, 
and change variables to $z=\sigma_u(x)$.
The inverse of the map $u\mapsto z$ and its Jacobian
on the tangent space of $\SSS^d$ are given by
$$
u(z;x) = \frac{x-z}{|x-z|}\,,\quad 
j(z;x)= \bigl(2|x-z|^{d-1}\bigr)^{-1} \,,
$$ 
where $|x-z|=2\sin\frac{\delta(x,z)}{2}$ is the Euclidean
distance between $x$ and $z$ in $\RRR^{d+1}$. We obtain
\begin{eqnarray*}
m\bigl( A\bigtriangleup A^*\bigr)-
\EE\bigl[ m\bigl(S_U A\bigtriangleup A^*\bigr)\bigr]
&=& 2 \int_{A^*\setminus A} \int_{A\setminus A^*} 
|x-z|^{-(d-1)}\, dm(z)dm(x)\\
&\ge& 2^{-d} \bigl(m( A\bigtriangleup A^*)\bigr)^2 \,.
\end{eqnarray*}
For a random sequence $S_{U_1\dots U_n}$, we 
take expectations again and apply Jensen's inequality to see that
\begin{eqnarray*}
\EE\Bigl[m(S_{U_1\dots U_{n-1}}A\bigtriangleup A^*)
- m(S_{U_1\dots U_n}A\bigtriangleup A^*)\Bigr]
&\ge&  2^{-d} \Bigl(\EE \bigl[m(S_{U_1\dots U_{n-1}}
A\bigtriangleup A^*)\bigr]\Bigr)^2\,.
\end{eqnarray*}
It follows that
$z_n=2^{-d} \EE\bigl[ m(S_{U_1\dots U_n}A\bigtriangleup  A^*)\bigr]$
satisfies $z_{n}^{-1}\ge z_{n-1}^{-1} + 1$,
proving Eq.~(\ref{eq:converge}) with constant $C_d=2^d$.
\hfill$\Box$

\section{Random compressions}

Let $A$ be a spherical cap centered at
a point $a$. Polarization with respect to a reflection
$\sigma:x\mapsto \bar x$ transforms $A$ into the 
spherical cap of the same volume centered at $\tau(a)$, where
$$
\tau(x) = \left\{\begin{array} {ll}
x \,,\quad &\delta(x,O)\le  \delta(\bar x,O)\,,\\
\bar x \,,\quad & \mbox{otherwise}\,.
\end{array}\right.
$$
We will refer to $\tau$ as the {\bf compression}
associated with $\sigma$.
The compression associated with a random reflection
$\sigma_U$ will be denoted by $\tau_U$.
The following lemma describes the 
distribution of the 
distance of $\tau_U(x)$ from the north pole.

\bigskip\noindent{\bf Lemma.} \ {\em 
If $U$ is uniformly distributed on $\Omega$,
then for every point $x\in\SSS^d$ with $\delta(x,O)=\xi$
\begin{equation} \label{eq:recursion}
P(\delta(\tau_U(x),O)> \beta)
= I_{\xi> \beta} \left\{
1-\frac{1}{\pi} \int_{0}^{\beta} 
\Bigl(\frac{\cos\theta-\cos\beta}{\cos\theta-\cos\xi}\Bigr)^{(d-1)/2}\,
d\theta\right\}\,,\qquad \beta\in [0,\pi]\,.
\end{equation}
}

\smallskip\noindent{\em Proof.} \ 
By definition of the compression, 
$$
P(\delta(\tau_U(x),O)> \beta)
= I_{\xi> \beta} \, P(\delta(\sigma_U(x),O)> \beta)\,.
$$
For $\xi\le \beta$, there is nothing more to show.
For $\xi>\beta$, we set $t=\cos\beta$
and calculate the spherical average as an expectation with 
respect to the standard normal probability measure 
on $\RRR^{d+1}$, see~\cite[Exercise 63 on p.80]{F}.
We use the coordinate system $u=(r\cos\theta, r\sin\theta, \hat u)\in 
\RRR\times\RRR\times \RRR^{d-1}$,
which we rotate into a position where
$x=\bigl(\cos\frac{\xi}{2},\sin\frac{\xi}{2},0\bigr)$
and $O=\bigl(\cos\frac{\xi}{2},-\sin\frac{\xi}{2}, 0\bigr)$.
Then
$$
(u\cdot x)(u\cdot O)= \frac{r^2}{2}(\cos 2\theta + \cos\xi)\,,
$$ 
and $\delta(\sigma_u(x),O)\le\beta$ if and only if 
$
-r^2(\cos 2\theta + \cos\xi)
\ge (|\hat u|^2+r^2)(t-\cos\xi)$.
This results in 
\begin{eqnarray*}
P(\delta(\tau_U(x),O)\le \beta)
&=&\int_{\SSS^d} I_{\sigma_u(x)\cdot O\ge t} \, dm(u) \\ 
&=& \frac{1}{2\pi}\int_{-\pi}^\pi
\int_{\RRR^{d-1}} \int_0^\infty   
I_{-r^2(\cos2\theta+t)\ge  |\hat u|^2(t-\cos\xi)}\,
2r e^{-r^2}\, dr\, d\gamma(\hat u)\,d\theta\,,
\end{eqnarray*}
where $\gamma$ is the standard normal probability measure 
on $\RRR^{d-1}$.
We integrate explicitly over
$r$ and evaluate the remaining Gaussian integral 
by rescaling $\hat v=\left(1-\frac{t-\cos\xi}{t+\cos 2\theta}\right)^{1/2}
\hat u$,
\begin{eqnarray*}
P(\delta(\tau_U(x),O)\le \beta)
&=& \frac{1}{2\pi}\int_{-\pi}^\pi 
I_{\cos2\theta +t< 0}
\left\{
\int_{\RRR^{d-1}} e^{\frac{t-\cos\xi}{t+\cos2\theta}|\hat u|^2}\, 
d\gamma(\hat u)\right\}\,d\theta\\
&=& \frac{1}{2\pi}\int_{-\pi}^\pi 
I_{\cos2\theta +t< 0}
\Bigl(\frac{-\cos2\theta-\cos\beta}{-\cos2\theta-\cos\xi}\Bigr)^{(d-1)/2}
\,d\theta\,.
\end{eqnarray*} 
The claim follows after restricting
the integral to a half-period 
and changing variables $2\theta\to\pi-\theta$.
\hfill $\Box$

\bigskip For $d=1$, 
the reflected point $\sigma_U(x)$ 
is uniformly distributed on $\SSS^1$, and
Eq.~(\ref{eq:recursion}) reduces to 
$$
P(\delta(\tau_U(x),O)> \beta)
= I_{\delta(x,O)> \beta} \left( 1-\tfrac{\beta}{\pi}\right)\,,
\qquad \beta\in [0,\pi]\,.
$$
As $d$ increases, $\sigma_U(x)$ concentrates
in a ball of radius comparable to $d^{-1/2}$ about $x$, 
its distance from the north pole concentrates in an interval 
of length comparable to $d^{-1}$ about $\xi$, 
and the integral in Eq.~(\ref{eq:recursion}) goes to zero.  
For all $d\ge 1$ and $0\le \beta\le \xi\le \pi$, we have the bound
$$
\frac{1}{\pi}\int_0^\beta \Bigl(\frac{\cos\theta-\cos\beta}
{\cos\theta-\cos\xi}\Bigr)^{(d-1)/2} d\theta 
\ \le\  \frac{\beta}{\pi} 
\Bigl(\frac{1-\cos\beta}{1-\cos\xi}\Bigr)^{(d-1)/2}
= \ (1+{\cal O}(\beta^2))\frac{2\sin\tfrac{\beta}{2}}{\pi}
\biggl(\frac{\sin\frac{\beta}{2}}{
\sin\frac{\xi}{2}}{ }\biggr)^{d-1}\,.
$$
If $\{U_i\}_{i\ge 1}$ is a sequence of independent uniformly
distributed random variables in $\Omega$,
it follows that the
Euclidean distance $Y_n=|\tau_{U_1\dots U_n}(x)-O|$
satisfies the recursion
\begin{equation} \label{eq:order}
P(Y_{n+1}>\eta \,\vert\, Y_n)\ge
I_{Y_n>\eta} \Bigl\{1-\frac{\eta}{\ell}\Bigl(\frac{\eta}{Y_n}\Bigr)^{d-1}
\Bigr\}\,,\qquad \eta\in [0,\ell]
\end{equation}
with initial value $Y_0=|x-O|=2\sin\frac{\xi}{2}$
and with $\ell=\pi-{\cal O}(\xi^2)$.

\section{The lower bound on the sharp constant}

Let $A$ be the hemisphere centered at a point
$a\ne O$, and set $\alpha=\delta(a,O)$.
We claim that
\begin{equation}\label{eq:liminf}
\liminf_{n\to\infty}~n\, \EE\bigl[m(S_{U_1\dots U_n}A\bigtriangleup A^*) \bigr] 
\ge (1-{\cal O}(\alpha^2))\,d\,.
\end{equation}
Taking $\alpha\to 0$, we see that
the sharp constant satisfies
$C_d\ge d$, completing the proof of the proposition.

To prove the claim,
consider a sequence of random points $\{V_i\}_{i\ge 1}$
that are distributed independently and uniformly 
on an interval $[0,\ell]$, and let $\tilde Y_n$ be the
$d$-th lowest point among $V_1,\dots, V_{n+d}$.
The random variable $\tilde Y_n$ is
called the {\bf $d$-th order statistic}
of $V_1\dots, V_{n+d}$.  
The sequence $\{\tilde Y_n\}_{n\ge 0}$ 
solves Eq.~(\ref{eq:order}) with equality,
because conditioned on $\tilde Y_n=y$, 
the $d\!-\!1$ points among $V_1,\dots V_{n+d}$ 
to the left of $y$ are independent and uniformly 
distributed on $[0,y]$, and $V_{n+d+1}$ is independent 
and uniformly distributed on $[0,\ell]$. The joint
distribution of the order statistics 
can be written explicitly in terms of binomial random 
variables $B(n,p)$, see~\cite[Exercises 21-25 on p.~142]{GS}.
We have
$$
P(\tilde Y_n>\eta\,\vert\, \tilde Y_0=y) =
I_{y>\eta}\sum_{j+k<d} P\bigl(B(d\!-\!1,\tfrac{\eta}{y})=j\bigr)
\cdot P\bigl(B(n,\tfrac{\eta}{\ell})=k\bigr)\,,
$$
where the first factor in  the sum accounts for the
points among $V_1,\dots V_d$
that fall to the left of $\eta$, while
the second factor accounts for 
such points among $V_{d+1},\dots V_{n+d}$.
By Stirling's formula,
$$
P(n\tilde Y_n>\eta \,\vert\, \tilde Y_0=y)\ \to\  
P\left(\Gamma(d)>\tfrac{\eta}{\ell}\right)\qquad (n\to\infty)
$$
for each $y\in (0,\ell]$, where $\Gamma(d)$ 
is a Gamma random variable that
describes the $d$-th point in a Poisson process 
of intensity one~\cite[Exercise 24 (b) on p.142]{GS}. 
In particular, $\EE[\tilde Y_n\,\vert\,\tilde Y_0=y]\to \ell d$.

The center of $S_{U_1\dots U_n}A$ is given by $\tau_{U_1\dots U_n}(a)$.
We have shown in Section~3 that $Y_n=|\tau_{U_1\dots U_n}(a)-O|$ satisfies
Eq.~(\ref{eq:order}). Since the right hand side of this recursion 
increases with $Y_n$ and the geodesic distance on the sphere 
exceeds the Euclidean distance, 
$$
P(\delta(\tau_{U_1\dots U_n}(a),O)>\eta) 
\ge P(\tilde Y_n>\eta\,\vert\, \tilde Y_0=\alpha)\,
$$
for all $n\ge 0$ with $\ell=\pi-{\cal O}(\alpha^2)$. 
For the mean, this implies that
$$
\liminf_{n\to\infty}
n\, \EE\bigl[\delta(\tau_{U_1\dots U_n (a)}, O)\bigr]
\ge (\pi-{\cal O}(\alpha^2))\,d\,.
$$
Eq.~(\ref{eq:liminf}) follows because 
the symmetric difference between two hemispheres
is just the distance of their centers, expressed
as a fraction of $\pi$.
\hfill$\Box$

\bigskip\noindent{\bf Remark.} \ 
A slightly more careful analysis of Eq.~(\ref{eq:recursion})
shows that for $a\ne O$,
$$
n \, \delta(\tau_{U_1\dots U_n}(a),O)
\to \pi \,\Gamma(d)\qquad (n\to\infty)
$$
in distribution,
and hence $\displaystyle{\lim_{n\to\infty}}~
n\,\EE[m(S_{U_1\dots U_n}A\bigtriangleup A^*)]=d$ 
for the hemisphere $A$ centered at $a$.


\end{document}